\documentclass[12pt]{article}%
\usepackage{amsfonts}
\usepackage{amsmath,amscd,amsthm,a4,amssymb}
\usepackage{amsmath,amscd,amsthm,a4,amssymb}
\usepackage{amsmath}
\usepackage{amssymb}
\usepackage{graphicx,multirow,rotating}%

\providecommand{\U}[1]{\protect\rule{.1in}{.1in}}

{}

{}

{}

{}

{}
\theoremstyle{plain}
{}

{}

{}

{}

{}

{}

{}

{}

{}

{}

{}

{}

{}

{}

{}

{}
\begin{document}

\begin{center}
{\Large \textbf{A New Descent Algebra for Standard Parabolic Subgroups of $W(A_n)$}}%

\end{center}

\centerline{\large T\"{u}lay YA\v{G}MUR $^{1}$, Himmet CAN $^{2}${\footnotetext{
{E-mail: $^{1}$tyagmur@erciyes.edu.tr (T.Ya\v{g}mur); $^{2}$can@erciyes.edu.tr (H.Can)}} }}

\ 

\centerline{\it $^{1,2}$ Department of Mathematics, Erciyes University, {\textsc{38039}}
Kayseri, Turkey}

\begin{abstract}
A new descent algebra $\sum_{W}(A_{n})$  of Weyl groups of type $A_n$, constructed by present authors in [1], is generated by equivalence classes $[x_J]$ arising from the equivalence relation defined on the set of all $x_J$. In this  paper, we introduce the structure of this new descent algebra for standard parabolic subgroups $W_J$ of Weyl groups of type $A_n$.
\end{abstract}

\begin{quote}
{\small \textit{{Key words and Phrases}: Weyl groups, descent algebra.} }

{\small 2010 \textit{Mathematics Subject Classification}: 20F55, 20F99. }
\end{quote}
 \bigskip
\section{Introduction}

\indent In [1], the authors  have constructed a new descent algebra of Weyl groups of type $A_n$.  This new descent algebra is denoted   by $\sum_{W}(A_{n})$. The main objective of this paper is to introduce and study the structure of this new descent algebra when the group is taken a standard parabolic subgroup of Weyl group of type $A_n$. \\
\indent In the next section, we shall give the structure of  a new descent algebra for standard paraolic subgroups $W_J$ of $W(A_n)$. This algebra will be denoted by  $\sum_{{W_J}}(A_{n})$. Additionally, we shall define the  induction and restriction maps between new descent algebras $\sum_{{W_J}}(A_{n})$ and $\sum_{{W_M}}(A_{n})$, where $W_J$ and $W_M$ are standard parabolic subgroups of $W(A_n)$.\\

\indent We now give the notation, which is fairly standard and follows that given in Ya\v{g}mur and Can [1], Carter [2], Solomon [3] and Bergeron et. al [4].\\
\indent Let  $e_{1}, e_{2}, ... , e_{n+1}$   be an orthonormal basis of a Euclidean space of dimension $n+1$.
Then the simple system and root system of type $A_{n}$ are given by 
\begin{center}
 $\Pi =\{e_{1}-e_{2}, e_{2}-e_{3}, ...  , e_{n-1}-e_{n}, e_{n}-e_{n+1}\}$, 
\end{center}

\begin{center}
 $\Phi =\{e_{i}-e_{j} | {i}\neq{j}, \quad  i,j=1,...,n+1 \}$,
\end{center}
respectively.\\
\indent  Let $\Pi$ be a simple system in root system $\Phi$ of type $A_{n}$ and $\Phi^{+}$ be the 
corresponding positive system. Then $W(A_{n})$ is called the {\textit{ Weyl group of type
$A_{n}$}} generated by the reflections $w_{r}$ for all $r\in\Phi$ [2]. For simplicity, throughout this study, this group will be denoted by $W$. \newline
\indent  Let $W_{J}$ be the subgroup of $W$, where $J$ is any subset of $\Pi$.
This group is called a standard parabolic subgroup of $W$.  Let $X_{J}$  be the
set of  representatives of the  cosets of $wW_{J}$ in $W$. 
Then $X_{JK}=X_{J}^{-1}\bigcap{X_{K}}$ is a set of distinguished double coset representatives of $W_{J}wW_{K}$  in $W$.\newline
\indent Let $x_J=\sum_{d\in{X_J}}{d}$. The set of all $x_J$ is a basis for an algebra $\sum(W)$ over the field of rationals with integer structure constants  $a_{JKL}$, where $a_{JKL}$  is the number of elements $w\in{X_{JK}}$ such that $w^{-1}(J)\cap{K}=L$. This algebra, discovered by Solomon [3] in 1976, is called {\textit{ the descent algebra (or Solomon algebra) of Weyl groups $W$}}. Moreover, in 1992 Bergeron et al. [4] has reconstructed the descent algebra systematically. \newline
\indent In [1], present authors have defined an equivalence relation on the set of all $x_J$ in order to form a basis for a new descent algebra  $\sum_{W}(A_{n})$. Furthermore, the basis of this new descent algebra consists of equivalence classes $[x_J]$ arising from the equivalence relation 
\begin{eqnarray}
x_{J}\sim{x_{K}}\Leftrightarrow{J\sim{K}}\Leftrightarrow{K=w(J),\quad
for\quad w\in{W}}.\nonumber
\end{eqnarray}
Additionally, a ring multiplication for two basis elements is defined by 
  \begin{eqnarray}
 [x_{J}][x_{K}]=[x_{J}x_{K}]=\sum_{L\subseteq{K}}a_{JKL}[x_{L}],
\end{eqnarray}
where $a_{JKL}$'s are defined as Solomon' descent algebra.\\

\indent The following  theorem is proved in [1, Theorem A].\newline

\noindent\textbf{Theorem 1.1} {\textit{Let $J\subseteq{\Pi}$, then
\begin{center}
$\sum_{W}(A_{n})=Sp\{ [x_J] | J\subseteq{\Pi}\}$
\end{center}
 is a new descent algebra of Weyl groups of type $A_{n}$.  }}\\

\section{Main Results}
 In this section we firstly define a new descent algebra for standard parabolic subgroups $W_J$ of Weyl groups of type $A_n$. After,  we shall study on this new algebra.  \\
\indent Let $K\subseteq{J}\subseteq{\Pi}$. An element $[x_{K}^J]$ can be defined as
\begin{eqnarray}
[x_{K}^J]=\{x_{L}^J|\quad L\sim{K},\quad L\subseteq{J}\}.\nonumber
\end{eqnarray}
\ Note that, if we take $J=\Pi$, then  because of $x_{L}^{\Pi}=x_L$ (see Bergeron et. al [4]) we obtain 
\begin{align}
[x_{K}^{\Pi}]  & =\{x_{L}^{\Pi}|\quad L\sim{K},\quad L\subseteq{\Pi}\}\nonumber\\
&=\{x_{L}|\quad L\sim{K}\}\nonumber\\
&=[x_K].\nonumber
\end{align}

For this reason, 
\begin{center}
$\sum_{{W_J}}(A_{n})=Sp\{ [x_{K}^J] | K\subseteq{J}\}$
\end{center}
is a new descent algebra for parabolic subgroups of Weyl groups of type $A_n$. Furthermore, the set $\{ [x_{K}^J] | K\subseteq{J}\}$ is a basis for this algebra.\\

\noindent \textbf{Theorem 2.1} {\textit{ Let $a_{MNP}^J$  be the structure constants of the new descent algebra                    $\sum_{{W_J}}(A_{n}) $  corresponding to the $[x_{K}^J]$ basis.
If $K,N\subseteq{\Pi}$, then 
\begin{eqnarray}
[x_K][x_N]=\sum_{P\subseteq{N}}{(\sum_{M\subseteq{J}}{a_{KJM}}{a_{MNP}^J})}[x_P],\nonumber
\end{eqnarray}\\
for all $J\subseteq{\Pi}$ such that $N\subseteq{J}$. Moreover, the structure constants satisfy the identities
\begin{eqnarray}
a_{KNP}=\sum_{M\subseteq{J}}{a_{KJM}}{a_{MNP}^J}.\nonumber
\end{eqnarray}
}}

For the proof of this theorem, we need to state the following crucial lemma proved in [4, Lemma 2.1].\\

\noindent \textbf{Lemma 2.2} {\textit{ If $K\subseteq{J}\subseteq{\Pi}$, then $x_K={x_J}{x_{K}^J}$.  
}}\\
\bigskip

\noindent \textbf{Proof of Theorem 2.1}
For $K,N\subseteq{\Pi}$ and $N \subseteq{J}$, by using Lemma 2.2 and equation (1), we have
\begin{align}
[x_K][x_N] & =[{x_K}{x_N}]\nonumber\\
&=[{x_K}{x_J}{x_{N}^J}]\nonumber\\
&=[{x_K}{x_J}][{x_{N}^J}]\nonumber\\
&=\sum_{M\subseteq{J}}{a_{KJM}}{[x_M][{x_{N}^J}]}\nonumber\\
&=\sum_{M\subseteq{J}}{a_{KJM}}{[x_J{x_{M}^J}][{x_{N}^J}]}\nonumber\\
&=\sum_{M\subseteq{J}}{a_{KJM}}{[x_J][{x_{M}^J}{x_{N}^J}]}\nonumber
\end{align}
\begin{align}
&=\sum_{M\subseteq{J}}{a_{KJM}}{[x_J]}{\sum_{P\subseteq{N}}{a_{MNP}^J}{[{x_{P}^J}]}}\nonumber\\
&=\sum_{{M\subseteq{J}},{P\subseteq{N}}}{a_{KJM}}{a_{MNP}^J}{[x_J]}{[{x_{P}^J}]}\nonumber\\
&=\sum_{{M\subseteq{J}},{P\subseteq{N}}}{a_{KJM}}{a_{MNP}^J}{[x_J{x_{P}^J}]}\nonumber\\
&=\sum_{{M\subseteq{J}},{P\subseteq{N}}}{a_{KJM}}{a_{MNP}^J}{[x_P]}\nonumber\\
&=\sum_{P\subseteq{N}}{(\sum_{M\subseteq{J}}{a_{KJM}}{a_{MNP}^J})}[x_P].\nonumber
\end{align}
Moreover, by the equation (1), we get 
\begin{eqnarray}
[x_K][x_N]=\sum_{P\subseteq{N}}{a_{KNP}}[x_P]\nonumber
\end{eqnarray}\\
So, since the set $\{[x_P]| P\subseteq{N}\}$ is linear independence, we obtain
\begin{eqnarray}
a_{KNP}=\sum_{M\subseteq{J}}{a_{KJM}}{a_{MNP}^J}.\nonumber
\end{eqnarray}\\
This completes the proof. \quad \quad\quad \quad \quad\quad \quad  \quad \quad \quad\quad \quad \quad \quad \quad \quad \quad \quad \quad \quad \quad \quad \quad  $\Box$
\bigskip

 We now give the definitions of  induction and restriction maps between the new descent algebras $\sum_{{W_J}}(A_{n})$ and $\sum_{{W_M}}(A_{n})$, but for simple presentation, we may denote these algebras as $\sum_{{W_J}}$ and $\sum_{{W_M}}$, respectively.\\

\noindent \textbf{Definition 2.3} Let $J\subseteq{M}\subseteq{\Pi}$. \textit{The induction map} between $\sum_{{W_J}}$ and $\sum_{{W_M}}$ is defined by
\begin{center}
${ind}_J^M: \sum_{{W_J}}\longrightarrow{\sum_{{W_M}}}$ ; \quad ${{ind}_J^M([{x_{K}^J}])=[{x_{J}^M}][{x_{K}^J}]}$.
\end{center}

\noindent \textbf{Definition 2.4} Let $J\subseteq{M}\subseteq{\Pi}$. \textit{The restriction map} between $\sum_{{W_M}}$ and $\sum_{{W_J}}$ is defined by
\begin{center}
${res}_J^M: \sum_{{W_M}}\longrightarrow{\sum_{{W_J}}}$ ;  
\end{center}
\begin{eqnarray}
{{res}_J^M([{x_{K}^M}])=\sum_{d\in{{W_M}\bigcap{X_{KJ}}}}[{x_{J\cap{d^{-1}(K)}}^J}]=\sum_{d\in{{W_M}\bigcap{X_{JK}}}}[{x_{J\cap{d(K)}}^J}]}.\nonumber
\end{eqnarray}

In [1, Proposition 2.10 ], the authors have proved that there is an isomorphism between  $\sum_{W}(A_{n})$ and $PB(W)$, where $PB(W)$ is the parabolic Burnside ring of associated Weyl group. For the case of induction, the permutation representation $W_{J}/W_{K}$ in $PB(W_J)$ induced to $PB(W)$ is simply $W/W_{K}$ [4].
The following theorem is arising from this discussion.\\

\noindent \textbf{Theorem 2.5} {\textit{ Given $J\subseteq{M}\subseteq{\Pi}$, let 
\begin{center}
$\Theta_J:\sum_{{W_J}}\to{\sum_{{W_M}}}$  and $\Theta_M:PB({W_J})\to{PB({W_M})}$  
\end{center}
be the canonical homomorphisms. Then,
\begin{center}
${\Theta_J}{res_J^{M}}={res_{W_{J}}^{W_{M}}}{\Theta_M}$\\
\end{center}
and
\begin{center}
${\Theta_M}{ind_J^{M}}={ind_{W_{J}}^{W_{M}}}{\Theta_J}$,\\
\end{center}
where
\begin{center}
${res_{W_{J}}^{W_{M}}}: PB({W_M})\to{PB({W_J})}$ ;
\end{center}
\begin{eqnarray}
{res_{W_{J}}^{W_{M}}}(W_{M}/W_{K})=\sum_{d\in{{W_M}\bigcap{X_{KJ}}}}W_{J}/W_{J\cap{d^{-1}(K)}}\nonumber
\end{eqnarray}
and
\begin{center}
${ind_{W_{J}}^{W_{M}}}: PB({W_J})\to{PB({W_M})}$ ; ${ind_{W_{J}}^{W_{M}}}(W_{J}/W_{K})=W_{M}/W_{K}$.
\end{center}
}}
\bigskip
\noindent \textit{Proof}. For all $[x_K^{M}]\in{\sum_{{W_M}}}$, by using Definition 2.4 and considering the existence of isomorphism between $\sum_{{W_M}}$ and $PB({W_M})$ (see Ya\v{g}mur and Can [1]), we get
\begin{align}
{\Theta_J}{res_J^{M}}([x_K^{M}]) &= {\Theta_J}({res_J^{M}}([x_K^{M}]))\nonumber\\
&={\Theta_J}(\sum_{d\in{{W_M}\bigcap{X_{JK}}}}[{x_{J\cap{d(K)}}^J}])\nonumber\\
&=\sum_{d\in{{W_M}\bigcap{X_{JK}}}}{\Theta_J}([{x_{J\cap{d(K)}}^J}])\nonumber\\
&=\sum_{d\in{{W_M}\bigcap{X_{JK}}}}{W_{J}/W_{J\cap{d(K)}}}.\nonumber
\end{align}
On the other hand, we obtain

\begin{align}
{res_J^{M}}{\Theta_J}([x_K^{M}]) &={res_J^{M}}({\Theta_J}[x_K^{M}])\nonumber\\
&={res_J^{M}}(W_{M}/W_{K})\nonumber\\
&=\sum_{d\in{{W_M}\bigcap{X_{KJ}}}}W_{J}/W_{J\cap{d^{-1}(K)}}\nonumber\\
&=\sum_{d\in{{W_M}\bigcap{X_{JK}}}}{W_{J}/W_{J\cap{d(K)}}}.\nonumber
\end{align}
This proves the first assertion of theorem.\\
\indent Now, for all $[x_K^{J}]\in{\sum_{{W_J}}}$, by using Definition 2.3 and considering the existence of isomorphism between $\sum_{{W_J}}$ and $PB({W_J})$ (see Ya\v{g}mur and Can [1]), we have

\begin{align}
{\Theta_M}{ind_J^{M}}([x_K^{J}]) &= {\Theta_M}({ind_J^{M}}([x_K^{J}]))\nonumber\\
&={\Theta_M}([x_J^{M}][x_K^{J}])\nonumber\\
&={\Theta_M}([x_K^{M}])\nonumber\\
&=W_{M}/W_{K}.\nonumber
\end{align}
Besides, we obtain
\begin{align}
{ind_J^{M}}{\Theta_M}([x_K^{J}]) &= {ind_J^{M}}({\Theta_M}([x_K^{J}]))\nonumber\\
&= {ind_J^{M}}(W_{J}/W_{K})\nonumber\\
&=W_{M}/W_{K}.\nonumber
\end{align}
This completes the proof. \quad \quad\quad \quad \quad\quad \quad  \quad \quad \quad\quad \quad \quad \quad \quad \quad \quad \quad \quad \quad \quad \quad \quad  $\Box$\\

We would also like to point out that these  results are obtained from conformation of Bergeron et. al[4]'s study to new descent algebra for parabolic subgroups of $W(A_n)$.

\end{document}